\documentclass{amsart}

\usepackage{graphicx}
\usepackage{amssymb}
\usepackage{amsthm}
\usepackage{amsmath}
\usepackage{pstricks,pst-node}
\input xy
\xyoption{all}

\newtheorem{theorem}{Theorem}%[section]
\newtheorem{lemma}[theorem]{Lemma}

\newtheorem{corollary}[theorem]{Corollary}

\newtheorem{cnj}[theorem]{Conjecture}

\theoremstyle{definition}
\newtheorem{definition}[theorem]{Definition}
\newtheorem{remark}[theorem]{Remark}

\newcommand{\Z}{\mathbb{Z}}
\newcommand{\Q}{\mathbb{Q}}
\newcommand{\F}{\mathbb{F}}

\def\M{{\mathcal M}}
\def\G{{\mathcal G}}

\newcommand{\inv}{^{-1}}

             \def \cH {{\mathcal H}}

             \def \N {{\mathbb N}}

\def\F{\mathbb F}

\def\C{\mathbb C}
\def\Q{\mathbb Q}
\def\Z{\mathbb Z}

\def\inv{^{-1}}

\DeclareMathOperator{\Hom}{Hom}
\DeclareMathOperator{\Frob}{Frob}
\DeclareMathOperator{\Gal}{Gal}

\def\Ker{\mbox{Ker}}
\def\im{\mbox{Im}}

\newcommand{\GL}{{\rm GL}}

\newcommand{\SL}{{\rm SL}}

\begin{document}

%%%from my template
\title[Torsion in cohomology and Galois representations]{Torsion in the cohomology of congruence subgroups of $SL(4,\Z)$ and Galois representations}

\author{Avner Ash} \address{Boston College\\ Chestnut Hill, MA 02445}
\email{Avner.Ash@bc.edu} \author{ Paul E. Gunnells}
\address{University of Massachusetts Amherst\\ Amherst, MA 01003}
\email{gunnells@math.umass.edu} \author{Mark McConnell}
\address{Center for Communications Research\\ Princeton, New Jersey 08540}
\email{mwmccon@idaccr.org}

\thanks{AA wishes to thank the National Science Foundation for support
of this research through NSF grant DMS-0455240, and also the NSA
through grant H98230-09-1-0050.  This manuscript is submitted for
publication with the understanding that the United States government
is authorized to produce and distribute reprints.  PG wishes to thank
the National Science Foundation for support of this research through
NSF grant DMS-0801214.}



\maketitle

\begin{abstract}
We report on the computation of torsion in certain homology theories
of congruence subgroups of $\SL(4,\Z)$.  Among these are the usual
group cohomology, the Tate-Farrell cohomology, and the homology of the
sharbly complex.  All of these theories yield Hecke modules.  We
conjecture that the Hecke eigenclasses in these theories have attached
Galois representations.  The interpretation of our computations at the
torsion primes 2,3,5 is explained.  We provide evidence for our
conjecture in the 15 cases of odd torsion that we found in levels $\le
31$.
\end{abstract}

\section{Introduction}\label{intro}

In a series of papers \cite{AGM1, AGM2, AGM3} we have computed the
cohomology $H^{5}$ of certain congruence subgroups of $\SL(4,\Z)$
(with $5$ being the degree of most interest for reasons explained in
\cite{AGM1}).  The coefficients in the cohomology of these papers
consists of the trivial module $\C$, or its stand-in, $\F_p$ for a
large prime $p$.

We are now beginning a series of computations of the torsion in the
cohomology with $\Z$-coefficients, from which we can also deduce the
the cohomology with $\Z/p\Z $-coefficients for all primes $p$. 

We are primarily interested in the cohomology as a module for the
Hecke algebra and in the connections with Galois representations.  In
the future, we hope to deal with twisted mod $p$-coefficients and to
be able to test examples of the Serre-type conjectures enunciated in
\cite{Ash-Sinnott, ADPS, Herzig}.  Looking at the details of the
these conjectures convinces us that any non-Eisensteinian example
likely to be within the range of feasible computation will involve
non-trivial coefficients.

There are also conjectures that are a sort of converse of the
conjectures above.  The prototype of these is Conjecture B in
\cite{A1}.  We have been led, through our computational work, to
generalize these conjectures from the group cohomology itself to a
variety of related theories.

\begin{cnj}\label{genconj}
A Hecke eigenclass in any reasonable (co)homology theory of an
arithmetic subgroup of $\GL(m,\Z)$ should have a Galois representation
attached.
\end{cnj}

We will give a precise version, including what ``reasonable" means, in
Conjecture~\ref{conj1} of Section~\ref{conjs} below.

As of now, our programs compute with trivial coefficients.  As the
level $N$ of the congruence group gets large, the time required by the
computations also becomes large.  In this paper we will:

(1) report on our results for $H^5(\Gamma_0(N),\Z)$ over $\Z$ with
$N\leq 31$ and test Conjecture~\ref{conj1} where possible;

(2) clarify what it is we are actually computing when $p=2, 3, 5$, since
the group $\Gamma_0(N)$ contains torsion elements of orders $2, 3, 5$; and

(3) make some specific conjectures along the lines of
Conjecture~\ref{genconj} and relate them to each other.

Among our preliminary results in (1), we find torsion classes of
orders $2, 3, 5$.  For these, we are actually testing one of these new
conjectures, contained in Conjecture~\ref{conj1}.  (For reasons
explained below we cannot handle the prime $2$ at present.)

In fact, we will discuss these issues in the reverse order (3)--(2)--(1).

\section{Homology theories and Hecke operators}\label{homology}

In this section we rely heavily on Brown's book \cite{B}.  Our ``value
added" is showing that the Hecke algebra acts on the various kinds of
homology theories defined below, and that the main exact sequence
below is Hecke equivariant.
 
In \cite{B} the theory we discuss is developed for any virtual duality
group (see \cite[p.~229]{B} for the definition).  Any arithmetic group
is a virtual duality group, as proved by Borel and Serre \cite{B-S}.
So let $G$ be a reductive $\Q$-group, $\Gamma$ an arithmetic subgroup
of $G(\Q)$, and $S$ a subsemigroup of $G(\Q)$ such that $(\Gamma,S)$
is a Hecke pair.  Let $\cH = H(\Gamma, S)$ denote the Hecke algebra of
double cosets with $\Z$-coefficients, which we assume to be
commutative. The paper \cite{A1} contains an introduction to these
terms.

First recall that if $M$ is any $S$-module, $\cH$ acts naturally on
$H_\ast(\Gamma,M)$ and $H^\ast(\Gamma,M)$.  The action is given by
composing a twisted restriction map with a corestriction map.  See,
for example, \cite[Def.~1.7]{A1}.  

We let $St$ denote the dualizing module for $\Gamma$, also known as
the \emph{Steinberg module}.  It is isomorphic to $H^n(\Gamma,\Z\Gamma)$
where $n$ is the virtual cohomological dimension of $\Gamma$.

Choose a complete resolution $(F,P,\epsilon)$ for $\Gamma$.  Recall
\cite[p.~273]{B} this means that $F$ is an acyclic chain complex of
projective $\Z\Gamma$-modules together with an ordinary projective
resolution $\epsilon\colon P\to\Z$ over $\Z\Gamma$ such that $F$ and
$P$ coincide in sufficiently high dimensions.  By the results of Borel
and Serre already mentioned, a complete resolution for $\Gamma$
exists.

For any $\Gamma$-module $M$ set \cite[p.~277]{B}
$$
\widehat H^\ast(\Gamma,M) = H^\ast(\Hom_\Gamma(F,M)).
$$
(We use Brown's conventions for the definition of $\Hom$
\cite[p.~5]{B}.  Further $M$ is to be thought of as a complex
concentrated in dimension 0.)

Up to isomorphism $H^\ast(\Hom_\Gamma(F,M))$ is independent of the
choice of $F$, and it is called the \emph{Farrell cohomology} of $\Gamma$.
Since restriction and corestriction maps exist in the theory of
Farrell cohomology, the action of $\cH$ on it is defined in the usual
way.

Besides the ordinary cohomology and the Farrell cohomology, we need to
consider a third homology theory, which we will refer to as \emph{Steinberg
homology.}  This is defined by \cite[p.~279]{B}
$$
\widetilde H_\ast(\Gamma,M) = H_\ast(\Gamma, St\otimes_\Z M)).
$$
Again, $\cH$ acts on this in the usual way using restriction
and corestriction maps.

These three homology theories fit together in a long exact sequence.
To see that this sequence is $\cH$-equivariant, we must recall how it
is obtained \cite[p.~280]{B}.  We choose a finite type projective
resolution $\epsilon\colon P\to \Z$ and a finite type projective resolution
$\eta\colon Q\to St$, both over $\Z\Gamma$.  If $A$ is any complex of
$\Z\Gamma$-modules, denote by $\overline A =
\Hom_{\Z\Gamma}(A,\Z\Gamma)$ the dual complex.  We let $\Sigma^k A$
denote the $k$-th suspension of $A$.

Recall that $n$ denotes the virtual cohomological dimension of
$\Gamma$.  Proposition (2.5) of \cite{B} states that we may choose the
complete resolution $F$ so that $\overline F$ is the mapping cone of a
chain map $\phi\colon \Sigma^{-n}Q \to \overline P$.

It follows immediately as on \cite[p.~280]{B}, that for any coefficient
module $M$, the exact homology sequence of this mapping cone yields
the long exact sequence
$$
\cdots \to \widetilde H_{n-i} \to H^i \to \widehat H^i \to \widetilde H_{n-1-i} \to H^{i+1} \to \cdots
$$
\begin{theorem}\label{bs}  
Assume that $St$ has a resolution by $S$-modules that are projective
as $\Gamma$-modules.  Then this exact sequence is equivariant for the
action of the Hecke algebra $\cH$.
\end{theorem}

\begin{proof}
Up to isomorphism, this exact sequence doesn't depend on the choice
of resolutions $P$ and $Q$. Unfortunately, \cite{B} assumes that
$P$ and $Q$ are of finite type, i.e. in each degree they are finitely
generated as $\Z\Gamma$-modules.  We will need to use projective
resolutions $\epsilon'\colon P'\to \Z$ and $\eta'\colon Q'\to St$ over $\Z\Gamma$
that are not of finite type.  So we first want to show that we can
replace $P$ by $P'$ and $Q$ by $Q'$ in the construction of the mapping
cone and yet derive the same exact sequence.

To see this, first keep $P$ fixed.  Given $Q$ and $Q'$ as above,
choose a homotopy equivalence $f\colon Q'\to Q$.  Define the chain map
$\phi'\colon \Sigma^{-n}Q' \to \overline P$ by $\phi'=\phi\circ f$.  Let the
mapping cone of the chain map $\phi\otimes_\Gamma
M\colon \Sigma^{-n}Q\otimes_\Gamma M \to \overline P\otimes_\Gamma M$ be
denoted by $C$.  Let the mapping cone of the chain map
$\phi'\otimes_\Gamma M\colon\Sigma^{-n}Q'\otimes_\Gamma M \to \overline
P\otimes_\Gamma M$ be denoted by $C'$.  We have the following diagram
of complexes, where the horizontal lines are exact and the vertical
arrow $g$ is induced by $f$:
$$
\xymatrix{
 0 \ar[r] & \overline P\otimes_\Gamma M\ar[r] & C \ar[r] &  \Sigma\Sigma^{-n}Q\otimes_\Gamma M  \ar[r] & 0\\
 0 \ar[r] & \overline P \otimes_\Gamma M \ar@{=}[u] \ar[r] & C' \ar[u]_g \ar[r] &  \Sigma\Sigma^{-n}Q'\otimes_\Gamma M   \ar[u]_{\Sigma\Sigma^{-n}f\otimes_\Gamma 1} \ar[r] & 0
}
$$
Taking homology we obtain two long exact sequences as follows: using
the fact that $f$ and hence ${\Sigma\Sigma^{-n}f\otimes_\Gamma 1}$ are
homotopy equivalences:
$$
\xymatrix{
\cdots {}\ar[r]& H^i \ar[r] & \widehat H^i \ar[r] &  
\widetilde H_{n-1-i}\ar[r]&\cdots\\
\cdots {}\ar[r]& H^i \ar@{=}[u] \ar[r] & H_i(C') \ar[u]_{g_\ast} \ar[r] &  \widetilde H_{n-1-i}  \ar[u]_{h} \ar[r]&\cdots
}
$$
Since $f$ and hence ${\Sigma\Sigma^{-n}f\otimes_\Gamma 1}$ are
homotopy equivalences, the maps marked $h$ are isomorphisms.  Hence so
are the maps marked $g_*$.

Therefore $H_i(C') \approx \widehat H^i$ and we can compute Farrell
cohomology in this way, using $Q'$ in place of $Q$. We now keep $Q'$
fixed and repeat the argument with $P'$ in place of $P$, using the
fact that if $x\colon P\to P'$ is a homotopy equivalence then
$\overline x\colon\overline P'\to\overline P$ is a homotopy
equivalence.  We deduce that we can compute Farrell cohomology and the
long exact sequence from the mapping cone, using $P',Q'$ in place of
$P,Q$ respectively.

Now we choose $P'$ and $Q'$ so they are complexes of $S$-modules,
where the $S$-module structure extends that of $\Gamma$.  For example,
we can let $P'$ be the standard resolution of $\Z$ for the group
$G(\Q)$.  The existence of $Q'$ is given by hypothesis.

We now have the exact sequence of $S$-module complexes:
$$
\xymatrix{
0 \ar[r] & \overline P' \otimes_\Gamma M \ar@{=} \ar[r] & C'  \ar[r] &  \Sigma\Sigma^{-n}Q'\otimes_\Gamma M    \ar[r] & 0
}
$$
If $X$ is an $S$-module, or a complex of $S$-modules, then
$X\otimes_\Gamma M = H_0(\Gamma, X\otimes_\Z M)$ has a natural
$\cH$-action, since it is a homology group of $\Gamma$ with $S$-module
coefficients.  Since $\cH$ thus acts on $\overline P' \otimes_\Gamma
M$ and $\Sigma^{-n}Q'\otimes_\Gamma M$, it also acts on
$\Sigma\Sigma^{-n}Q'\otimes_\Gamma M$ and on the mapping cone $C'$.
Taking homology of the complexes in the last exact sequence, one
obtains in this way a $\cH$-action on the homology groups which
coincides with the action on these groups as already defined.
Following out the proof of the snake lemma one easily sees that all
the maps in the long exact sequence of homology are $\cH$-equivariant.
\end{proof}

We have not attempted to prove in general that $St$ possesses a
resolution $Q'$ by projective $\Gamma$-modules that are also
$S$-modules.  However, in the case when $G=\GL(m)/\Q$, we can
construct such a $Q'$ as follows: let $Sh_{\bullet}$ be the sharbly complex,
whose definition is recalled in Section~\ref{what} below.  It is a
resolution of $St$ by $S$-modules, but it is not projective if
$\Gamma$ has nontrivial torsion.  However, the tensor product of any
$\Gamma$-module $M$ with $\Z\Gamma$ becomes a free $\Z\Gamma$-module,
as long as $M$ is $\Z$-free \cite[Cor.~5.7]{B}.  So we can take the
tensor product of $Sh_{\bullet}$ with a free resolution $P'$ of $\Z$ by
$S$-modules: $Q'=P'\otimes_\Z Sh_{\bullet}$.  This proves the following:

\begin{corollary}\label{cor-bs}  
Let $(\Gamma,S)$ be a Hecke pair contained in $\GL(m,\Q)$.
Then for any $S$-module $M$, the exact sequence
$$
\cdots \to H_{n-i}(\Gamma,St\otimes_\Z M) \to H^i(\Gamma,M) \to
\widehat H^i(\Gamma,M) \to H_{n-1-i}(\Gamma,St\otimes_\Z M) \to \cdots
$$
is equivariant for the action of the Hecke algebra $\cH$.
\end{corollary}

The last homology theory we consider is \emph{sharbly homology}
$H_\ast(\Gamma,Sh_{\bullet}\otimes_\Z M)$.  Its relationship to
$H_\ast(\Gamma,St\otimes_\Z M)$ will be considered in
Section~\ref{what}.  It is also naturally a $\cH$ module, since
$Sh_{\bullet}$ and $M$ are $S$-modules.

As we have already noted, if $\Gamma$ possesses nontrivial torsion,
then $Sh_{\bullet}$ will not be a projective resolution., So when
$\Gamma$ is a subgroup of $\GL(4,\Z)$, we end up in actuality
computing $H_\ast(\Gamma,Sh_{\bullet}\otimes_\Z M)$, rather than
$H_\ast(\Gamma,St\otimes_\Z M)$, when $\Gamma$ possesses nontrivial
torsion.  We have seen that we could compute the Steinberg homology by
replacing $Sh_{\bullet}$ by $Q'=P'\otimes_\Z Sh_{\bullet}$, or by
rigidifying the sharbly complex in some other way.  However, any
method of doing this that we have considered increases the number of
cells in the relevant dimensions so much as to make actual computation
infeasible.

\section{Some conjectures and their interrelationships}\label{conjs}

In this section we state conjectures that state that Galois
representations are attached to Hecke eigenclasses in any of our four
homology theories.  We make a few remarks on the relationships between
these conjectures.  We state these conjectures for the congruence
subgroups $\Gamma_0(N)$ so we can be definite in our notation.  They
are easily modified for any Hecke pair $(\Gamma,S)$ in $\GL(m,\Q)$.

Let $\Gamma_0(N)$ be the subgroup of matrices in $\SL(m,\Z)$ whose
first row is congruent to $(*,0,\ldots,0)$ modulo $N$.  Define $S_N$
to be the subsemigroup of integral matrices in $\GL(m,\Q)$ satisfying
the same congruence condition and having positive determinant
relatively prime to $N$.

Let $\cH(N)$ denote the $\Z$-algebra of double cosets
$\Gamma_0(N)S_N\Gamma_0(N)$.  Then $\cH(N)$ is a commutative algebra
that acts on the cohomology and homology of $\Gamma_0(N)$ with
coefficients in any $\Z[S_N]$ module.  When a double coset is acting
on cohomology or homology, we call it a Hecke operator.  Clearly,
$\cH(N)$ contains all double cosets of the form
$\Gamma_0(N)D(\ell,k)\Gamma_0(N)$, where $\ell$ is a prime not
dividing $N$, $0\leq k\leq m$, and
$$D(\ell,k)=\left(\begin{matrix}
1&&&&&\cr&\ddots&&&&\cr&&1&&&\cr&&&\ell&&\cr&&&&\ddots&\cr&&&&&\ell\cr\end{matrix}\right)$$
is the diagonal matrix with the first $m-k$ diagonal entries equal to
1 and the last $k$ diagonal entries equal to $\ell$.  It is known that
these double cosets generate $\cH (N)$.  When we consider the double
coset generated by $D(\ell,k)$ as a Hecke operator, we call it
$T(\ell,k)$.

\begin{definition}\label{def:hp}
Let $A$ be a ring and $V$ an $\cH(N)\otimes_\Z A$-module. Suppose that
$v\in V$ is a simultaneous eigenvector for all $T(\ell,k)$ and that
$T(\ell,k)v=a(\ell,k)v$ with $a(\ell,k)\in A$ for all prime $\ell\not
| \ N$ and all $0\leq k\leq m$.  If
$$\rho\colon G_\Q\to \GL(m,A)$$ is a continuous representation of
$G_{\Q} = \Gal (\overline\Q/\Q)$ unramified outside $L\in \N$, and
\begin{equation}\label{eqn:hp}
\sum_{k=0}^{n}(-1)^k\ell^{k(k-1)/2}a(\ell,k)X^k=\det(I-\rho(\Frob_\ell)X)
\end{equation}
for all $\ell\not | \ LN$, then we say that $\rho$ is attached to $v$ (or
that $v$ corresponds to $\rho$). 
\end{definition}  

Let $p$ be a prime, and $\overline \F_p$ an algebraic closure of
$\F_p$.  Let $M$ be a $S_N$-module that is a finite-dimensional
vector space over $\F_p$ on which $S_N$ acts via its reduction modulo
$p$.  We call such a module an admissible $S_N$-module. Then we make
the following conjectures:

\begin{cnj}\label{conj1} 
Fix $p,m,N,M$ as above.  Let $v$ be a Hecke eigenclass in

(a) $H_{\ast}(\Gamma_0(N),St\otimes_\Z M)$, 

(b) $H^\ast(\Gamma_0(N),M)$,

(c) $\widehat H^\ast(\Gamma_0(N),M)$, or

(d)  $H_{\ast}(\Gamma_0(N),Sh_{\bullet}\otimes_\Z M)$.

Then there exists a continuous representation unramified outside $pN$
$$\rho\colon G_\Q\to \GL(m,\overline \F_p)$$
attached to $v$.
\end{cnj}

It follows immediately from the exact sequence of
Corollary~\ref{cor-bs} that if any two of Conjectures~\ref{conj1}
(a),(b),(c) hold, then the third one holds.  We can express this also
in the following way.

Let $\M(N)$ denote a set of representatives of all admissible
$\F_p[S_N]$-modules modulo isomorphism.  If $V$ is an
$\cH(N)\otimes_\Z \F_p$-module, let $\G(V)$ denote the set of all
Galois representations (modulo isomorphism) attached to Hecke
eigenvectors in $V$.

\begin{theorem}\label{dimshift} We have the following inclusions of
sets of Galois representations:

\begin{equation}\label{inclusion1}
\G(\oplus_{M\in\M(N)} \oplus_k \widehat H^k(\Gamma_0(N),M))
\subset
\G(\oplus_{M\in\M(N)} \oplus_k  H^k(\Gamma_0(N),M)).
\end{equation}

\begin{equation}\label{inclusion2}
\G(\oplus_{M\in\M(N)} \oplus_k H_k(\Gamma_0(N),St\otimes_\Z M))
\subset
\G(\oplus_{M\in\M(N)} \oplus_k  H^k(\Gamma_0(N),M)).
\end{equation}
\end{theorem}

\begin{proof}
As explained on p.~278 of \cite{B}, there is dimension shifting in
both directions on Farrell cohomology.  We wish to dimension shift
upwards.  To do this, we use induced modules.  Fix a torsionfree
subgroup $\Gamma'$ of finite index in $\Gamma_0(N)$.  Let $I(M)$
denote the induced module $\F_p\Gamma_0(N)\otimes_{\F_p\Gamma'} M$.
Let $K$ denote the kernel of the natural map $I(M)\to M$.  Then
$$
\widehat H^k(\Gamma_0(N),M))\approx \widehat H^{k+1}(\Gamma_0(N),K)).
$$
Now $I(M)$ is isomorphic to the module denoted
$I(\Gamma',\Gamma_0(N);M)$ in Definition 1.5 p.~240 of \cite{A1} and
in particular has the structure of $\F_p S_N$-module.  Therefore
$\widehat H^k(\Gamma_0(N),M))$ and $\widehat H^{k+1}(\Gamma_0(N),K))$
are isomorphic as $\cH$-modules.

Lemma 1.6 of \cite{A1} implies that $I(M)$ and $K$ are admissible.
Repeating this construction, we can dimension shift as high as we
like.  So if $n$ is the virtual cohomological dimension of
$\Gamma_0(N)$, we obtain
$$
\G(\oplus_{M\in\M(N)} \oplus_k \widehat H^k(\Gamma_0(N),M))
\subset
\G(\oplus_{M\in\M(N)} \oplus_{k>n}  \widehat H^k(\Gamma_0(N),M)).
$$
But above dimension $n$, the Farrell and ordinary cohomology are
isomorphic as $\cH$-modules.  This proves the inclusion in \eqref{inclusion1}.

The inclusion in \eqref{inclusion2} now follows from the exact sequence of
$\cH$-modules of Corollary~\ref{cor-bs}, as in the proof of Theorem
3.1 of \cite{A1}.
\end{proof}

Remark: the inclusion \eqref{inclusion1} is strict, in general.  For example, if
$m=p-1$, Theorem 0.2 of \cite{A1} plus Theorem 7.3 of \cite{B} implies
that all Galois representations attached to the Farrell cohomology of
$\Gamma_0(N)$ with any admissible coefficient module are of a very
simple type, namely, induced from a character.  On the other hand, if
now $m=4, p=5$ we can reduce modulo $5$ any of the Hecke eigenclasses
found in \cite{AGM1, AGM2, AGM3} that have non-induced Galois
representations attached to them.

We do not know if the inclusion \eqref{inclusion2}  is strict.  We guess that it is.

As explained in Section~\ref{what}, the conjecture we actually test
with our computations on $\GL(4)$ is Conjecture~\ref{conj1} (d) for
$m=4$.  In Section~\ref{what}, we will explain the relationship
between Conjecture~\ref{conj1} (d) and Conjecture~\ref{conj1} (a).

\begin{theorem}\label{c54}
In all cases, Conjecture~\ref{conj1} (b) implies
Conjecture~\ref{conj1} (a) and Conjecture~\ref{conj1} (c).  

Now suppose that $p,m,N$ are such that there is no $p$-torsion in
$\Gamma_0(N)$.

(1) Conjecture~\ref{conj1} (c) is trivially true since the Farrell cohomology vanishes.  

(2) Conjectures~\ref{conj1} (a) and (b) and (d) are equivalent.

\noindent Next suppose that $p=5,m=4$.
Then 

(3) Conjecture~\ref{conj1} (c) is true.  

(4) Conjectures~\ref{conj1} (a) and (b) are equivalent.
\end{theorem}

\begin{proof}
The first assertion follows immediately from Theorem~\ref{dimshift}.

Now suppose there is no $p$-torsion in $\Gamma_0(N)$.  Then the
Farrell cohomology vanishes identically \cite[Exercise 2, p.~280]{B}.
Also, by Lemma~\ref{shst} below, the Steinberg cohomology is
isomorphic to the ordinary group cohomology, and to the sharbly
cohomology, because $\Gamma_0(N)$ is $p$-torsionfree.  Because of the
long exact sequence of Corollary~\ref{cor-bs}, (1) implies (2).

Next suppose $p=5,m=4$.  Then Theorem 6.4.3 of \cite{A1} and
Theorem 7.3 of \cite{B} imply that there is a Galois
representation attached to any Hecke eigenclass in the Farrell
cohomology.  This implies (3), and (4) then follows from the long exact
sequence of Corollary~\ref{cor-bs}.
\end{proof}

If there is $p$-torsion in $\Gamma_0(N)$, then 
Lemma~\ref{shst} will give us a relationship between
Conjectures~\ref{conj1} (a) and (d).

\section{The primes $2, 3, 5, p>5$}\label{235}

The torsion primes for $\GL(4,\Z)$ are $2,3,5$.  Setting $n=4$ in our notation, we have:

\begin{lemma}\label{lemma-tors} 
For any $N \ge 1$, the subgroup $\Gamma_0(N)$ of $\SL(4,\Z)$ has $2$-
and $3$-torsion.  It has $5$-torsion if and only if $N$ is not
divisible by $25$ and every prime divisor $p$ of $N$ not equal to $5$
satisfies the condition $5 \mid  p-1$.
\end{lemma}

\begin{proof}
Since $\Gamma_0(N)$ contains subgroups isomorphic to $\GL(3,\Z)$, it
contains $2$- and $3$-torsion.  Now consider $5$-torsion.  Because the
class number of $\Q(\zeta_5)$ is $1$, any subgroup of $\SL(4,\Z)$ of
order $5$ is conjugate to the one generated by the element
$$
Z=\left(\begin{matrix}
0&0&0&-1\cr1&0&0&-1\cr0&1&0&-1\cr0&0&1&-1\cr\end{matrix}\right).
$$  
So $\Gamma_0(N)$ contains $5$-torsion if and only if there exists
$A\in \SL(4,\Z)$ such that $A\inv ZA \in \Gamma_0(N)$ iff there exists
a primitive vector $v\in \Z^4$ and an integer $\lambda$ prime to $N$
such that $vZ \equiv \lambda v$ modulo $N$.  Keeping in mind that $v$
must be primitive, one sees that this is so if and only if there
exists an integer $a$ (namely any $a$ congruent to $\lambda$ modulo
$N$) such that $1+a+a^2+a^3+a^4 = 0$.  The conclusion of the lemma
follows easily.
\end{proof}

We now describe our computations of $p$-torsion in terms of the
various~$p$.

The prime 2: Unfortunately, at present, the algorithm in
\cite{experimental} for reducing sharblies, which is essential for
computing the Hecke operators on the sharbly homology, involves
division by 2.  So we cannot compute the Hecke action on 2-torsion
classes at this time.  Since there is generally a lot of 2-torsion in
the homology, fixing the algorithm to remove this division by 2 is a
pressing desideratum that we plan to address in future work.

The prime 3: We have examples of 3-torsion that verify
Conjecture~\ref{conj1} (d).

The prime 5: We have examples of 5-torsion that verify
Conjecture~\ref{conj1} (d).  Whether or not $\Gamma_0(N)$ contains
$5$-torsion, Lemma~\ref{c54} shows that when $p=5$
Conjectures~\ref{conj1} (a) and \ref{conj1} (c) are equivalent.  When
$\Gamma_0(N)$ does not contain $5$-torsion, Theorem~\ref{c54} tells us
that Conjectures~\ref{conj1} (a) and \ref{conj1} (d) are
equivalent. If $\Gamma_0(N)$ does not contain $5$-torsion, the
situation is the same as in the next paragraph.

Primes $p > 5$: We have examples of $p$-torsion that verify
Conjecture~\ref{conj1} (d).  Since $\Gamma_0(N)$ does not contain
$p$-torsion, Conjectures~\ref{conj1} (a), (b) and (d) are all
equivalent by Lemma~\ref{c54}.  So in these cases we are truly
verifying the original conjecture B of \cite{A1}.

\section{What we are computing}\label{what}

Let $(\Gamma,S)$ be a Hecke pair contained in $\GL(m,\Q)$.

In \cite{AGM1} we explained in detail how to compute with the spectral sequence
$$
E_1^{p,q} = \sum_{\sigma\in\Sigma_p} H^q(\Gamma_\sigma,M_\sigma) \Rightarrow 
H^{p+q}_\Gamma(W,M).
$$
This spectral sequence was derived as in \cite{B} from the double complex 
$$
\Hom_\Gamma(P_\bullet, C^\bullet(W,M)).
$$
Here, we have chosen a resolution $P_\bullet \to \Z$ of $\Z$ by free
$\Z S$-modules, and $W$ is a contractible $\Gamma$-cell complex.  (In
practice it is the well-rounded retract \cite{ash80}.)

When the torsion in $\Gamma$ is invertible on $M$, then $E_1^{p,q}=0$
for $q>0$.  Then $E_1^{p,0}=H^{p}_\Gamma(W,M)$ is isomorphic to
$H^{p}(\Gamma,M)$. The Hecke algebra $\cH(\Gamma,S)$ acts on
$E_1^{p,0}$ naturally via this isomorphism.  This is what we computed
in \cite{AGM1, AGM2, AGM3}.

However, when the torsion in $\Gamma$ is not invertible on $M$, then
the higher rows of the spectral sequence do not vanish and the Hecke
algebra does not act on the individual terms of the spectral sequence,
because $S$ does not preserve the cellulation.  In other words,
$E_1^{p,0}$ is no longer computing the object in Conjecture
\ref{conj1}(b).  We will now show that $E_1^{p,0}$ is actually
computing the object in Conjecture \ref{conj1}(d), the sharbly
homology.  The papers \cite{AGM1, AGM2, AGM3} computed both
\ref{conj1}(b) and \ref{conj1}(d)---they were the same because we were
using $\C$ coefficients.  The present paper extends those papers,
working in the setting \ref{conj1}(d).

%For  background we begin by recalling the following theorem.

%\begin{theorem}\label{bs}  
%There is a natural map, given by cap product with a ``fundamental class" 
%$z \in H_\nu(\Gamma,\St)$:
%$H^j(\Gamma,M) \to H_{\nu-j}(\Gamma, \St\otimes M)$.  If the torsion in $\Gamma$ is invertible on $M$, it is an isomorphism.
%\end{theorem}
%
%\begin{proof}
%See Bieri-Eckmann and Borel-Serre, or Brown.
%\end{proof}

%I'M NOT SURE THERE IS A FUNDAMENTAL CLASS WHEN GAMMA IS NOT TORSION FREE.

Recall that the sharbly complex $Sh_\bullet \to St$ gives a resolution
of the Steinberg module by $\Z S$-modules.  But it is not
$\Gamma$-projective if $\Gamma$ is not torsionfree.

Facts about the sharbly complex may be found in \cite{A2}.  The term
$Sh_k$ is the $\Z \GL(m,\Q)$-module generated by symbols
$[v_1,\dots,v_{m+k}]$ where the $v_i$ are primitive vectors in $\Z^m$,
subject to the relations

(i) $[v_{\sigma 1},\dots,v_{\sigma
(m+k)}]=(-1)^\sigma[v_1,\dots,v_{m+k}]$ for all permutations $\sigma$;

(ii) $[v_1,\dots,v_{m+k}] = 0$ if $v_1,\dots,v_{m+k}$ do not span all
of $\Q^m$.

We can form the double complex $P_\bullet \otimes_\Gamma Sh_\bullet$.
From this we get in the usual way a first quadrant homology spectral
sequence:
$$
E^1_{p,q} = H_q(\Gamma,Sh_p\otimes_\Z M) \Rightarrow
H_{p+q}(\Gamma,St\otimes_\Z M).
$$

Let $\cH$ denote the Hecke algebra $\cH(\Gamma,S)$.
\begin{lemma}\label{shst}

(a) For any $1 \le r \le \infty$, $\cH$ acts on every term $E^r_{p,q}$
in the $E^r$ page of the spectral sequence and commutes with all
differentials.  The abutment morphism is equivariant for $\cH$.

(b) Let $d$ be the product of all the torsion primes of $\Gamma$.
Suppose that $d$ acts invertibly on $M$.  Then $E_1^{p,q}=0$ for
$q>0$.  Then the sharbly homology $H_p(Sh_p\otimes_{\Z\Gamma} M)$ is
isomorphic to the Steinberg homology $H_{p}(St\otimes_{\Z\Gamma}
M)$. This isomorphism is equivariant for $\cH$.
\end{lemma}

\begin{proof}
Since $P_\bullet$ is a resolution of $S$-modules, as is $Sh_\bullet$,
$S$ acts on every term in the $E^0$ page of the spectral sequence and
the differentials are $S$-module maps.  It follows that $\cH$ acts on
every term in the $E^1$ page of the spectral sequence and the
differentials are $\cH$-module maps. Then (a) follows immediately.

It follows directly from the definition of the sharbly complex that
for each $i$, $Sh_i$ as a $\Gamma$-module is isomorphic to a direct
sum of induced $\Gamma$-modules.  Indeed, let $R$ be a set of
representatives of the $\Gamma$-orbits on ``pure" $p$-sharblies,
i.e. on the set of symbols $[v_1,\dots,v_{m+p}]$ modulo the action by
permutation of the primitive vectors $v_j\in \Z^m$.  For $r\in R$, let
$\Gamma_r$ denote the stabilizer, which is a finite group.  Then
$Sh_p$ is isomorphic to $\oplus_{r\in R}\Z[\Gamma/\Gamma_r]$.

Then
$$
H_q(\Gamma,Sh_p\otimes_\Z M) \approx 
\oplus_{r\in R}H_q(\Gamma,\Z[\Gamma/\Gamma_r]\otimes_\Z M)
\approx
\oplus_{r\in R}H_q(\Gamma_r, M)
$$
the last isomorphism following by Shapiro's lemma.  Now assume $d$
acts invertibly on $M$. Since $d$ divides the order of $\Gamma_r$ for
every $r$, each term in the last direct sum is 0 if $q > 0$, and so
$H_q(\Gamma, Sh_p \otimes_{\Z} M) = 0$ if $q >0$.

Therefore, each term in the complex computing the sharbly
homology, namely $E_1^{p,0}=H_0(\Gamma,Sh_p\otimes_\Z M)$, is
isomorphic to the corresponding term of the complex that computes the
Steinberg homology, namely $H_{0}(\Gamma,St\otimes_\Z M)$.  The rest
of (b) follows from (a).
\end{proof}

\begin{remark}
When $d$ is not invertible on $M$, all we can assert is the
following: (1) $E^\infty_{p,0}$ is a sub-$\cH$-module of the sharbly
homology $H_p(Sh_{\bullet}\otimes_\Gamma M) = E^2_{p,0}$.  (2) $E^\infty_{p,0}$
is a quotient-$\cH$-module of the Steinberg homology
$H_p(St\otimes_\Gamma M)$.  These assertions follow from standard
facts about spectral sequences.
\end{remark}

Let us now consider what we actually computed in \cite{AGM1, AGM2,
AGM3}.  Let $W$ be the well-rounded retract for $\GL(4)$.  Consider a
cell $\sigma$ of dimension $d>0$ in $W$ with minimal vectors
$v_1,\dots,v_k$.  Dual to this cell is the Voronoi cell with the same
minimal vectors.  Then $k+d=10$ and $\sigma$ corresponds to the
sharbly $[v_1,\dots,v_k]$. Our computations involve only $d=4,5,6$.
(If $d=0$, there is a cell $\sigma$ with $k+d>10$.  For this cell, we
would have to use a simplicial subdivision of the dual Voronoi cell
and then convert to sharblies.  Although we have no need to compute in
this dimension, similar phenomena will appear widely for $\GL(m)$
with $m>4$.)

Let us call a sharbly $[v_1,\dots,v_k]$ such that $v_1,\dots,v_k$ are the
minimal vectors of a cell in $W$ a \emph{V-sharbly}.  Then in the range
$k=4,5,6$, the $\Z$-spans of these V-sharblies form a subcomplex of the
sharbly complex.  We will also call any element of this span a
V-sharbly.  Write $[\sigma]$ for the V-sharbly corresponding to the
cell $\sigma$.  When we view this sharbly inside the
$\Gamma$-coinvariants, as when computing the homology of $\Gamma$ in
$Sh_\bullet \otimes M$, we write it as $[\sigma]_\Gamma$.

We compute the bottom row of Brown's spectral sequence,
i.e.~$\Ker(d_1)/\im(d_1)$, at the $E_1^{5,0}$ node.  A typical
cohomology class is thus represented by a linear combination $\sum
a_\sigma \sigma$ where $a_\sigma\in M$ and $\sigma$ runs over a set of
representatives of $5$-chains modulo $\Gamma$.  We convert this to the
V-sharbly $\sum a_\sigma [\sigma]_\Gamma$.  We then compute the Hecke
operator $T=\Gamma s \Gamma = \sum s_\alpha \Gamma$ on this sharbly,
obtaining $\sum\sum a_\sigma [\sigma]_\Gamma|s_\alpha$.  The algorithm
of \cite{experimental} is then used to rewrite this last expression as
homologous to some V-sharbly, $\sum b_\sigma [\sigma]_\Gamma$.  This
is reconverted into a homology class in the homology of the first
spectral sequence at the same node, $\sum b_\sigma \sigma$ and we thus
get a matrix for $T$ in terms of the chosen basis of cycles.

When working away from the torsion primes of $\Gamma$, namely 2, 3 and
possibly 5 (cf.~Lemma \ref{lemma-tors}), each of the two spectral
sequences are zero above the first row, and they compute the same
thing, Hecke-equivariantly.  But at a torsion prime, there is no such
vanishing, and the Hecke operators don't act on Brown's spectral
sequence.

In effect, at primes that divide the torsion in $\Gamma$, we are
computing as Hecke-module, not the cohomology of $\Gamma$ with
coefficients in $M$, but rather the homology of $\Gamma$ with
coefficients in $Sh_{\bullet}\otimes M$ at the $(1,0)$ node. The connection
between this and the group cohomology, such as it is, may be seen by
Corollary~\ref{bs} and the remark following Lemma~\ref{shst}.

Since we have generalized (in Conjecture~\ref{conj1} (d)) the
conjecture that expects Galois representations to be attached to
Hecke-eigenclasses in the group cohomology, we can expect Galois
representations to be attached to the Hecke-eigenclasses we are
computing.  This indeed happens in the examples we have computed so
far---see Section~\ref{data}.

There are two possible problems with these computations, both of which
we conjecture do not arise in practice.  First, it may not be possible
to reduce $\sum\sum a_\sigma [\sigma]_\Gamma|s_\alpha$ to a homologous
cycle of V-sharblies.  However, if the algorithm of
\cite{experimental} always terminates (as it always has in practice)
such a reduction is always possible.

The second problem is that the sharbly cycle handed to us to perform
Hecke operations on, namely $\sum a_\sigma [\sigma]_\Gamma$, could be
a boundary in the sharbly complex even though it is not a boundary in
the $V$-part of the sharbly complex.  In this case, any Hecke
eigenvalues we compute would be spurious, since they would be the
``eigenvalues" of an operator on the 0-vector.  If we could show that
all elements in $Sh_2$ are homologous to $V$-sharblies (either by an
algorithm or some other way) then this problem is irrelevant.  We
conjecture that this problem doesn't happen, as suggested by our
data. The Hecke eigenvalues we compute never appear to be
nonsensical, and we are always able to attach Galois representations
to our putative Hecke eigenclasses.  If this problem were happening,
then we should be obtaining random numbers for supposed Hecke
eigenvalues.

\section{Computational data}\label{data}

In this section we present the results of our current torsion
computations.  This data comes in fact from computations of the
$(1,0)$ node of the spectral sequence computing the homology of
$Sh_{\bullet}\otimes M$, where the module $M$ is $\Z$ with $\Gamma$
acting trivially.  Recall that the conjectures above concern modules
that are $\F_{p}$-vector spaces.  By the universal coefficient theorem, the
same packages of Hecke eigenvalues we present also occur in the
homology with coefficients in the trivial module $\F_p$.

We computed homology and Hecke operators using the same techniques in
\cite{AGM1, AGM2, AGM3}.  In particular, as in \cite{AGM1} we use a
slight modification of the sharbly complex $Sh_{\bullet}$ from
Section~\ref{what} that includes the extra relation
\[
[v_{1},v_{2},\dotsc ,v_{m+k}]-[-v_{1},v_{2},\dotsc ,v_{m+k}] = 0.
\]
The resulting complex is homotopy equivalent to $Sh_{\bullet}$ as long
as $2$ is invertible in the coefficients.  This causes no trouble for
the current work, since as mentioned before we only report on
$p$-torsion for $p$ odd.

Table \ref{tab:dims} shows the levels $\leq 31$ that have $p$-torsion
for $p$ odd, and gives the dimension of the relevant homology group as
an $\F_{p}$-vector space.  We note that we have examples for
$p=3, 5, 7$, and that for $p=5$ we have examples both where $\Gamma$ has
$5$-torsion and where $\Gamma$ doesn't (Lemma \ref{lemma-tors}).  If a
level $N\leq 31$ doesn't appear in Table \ref{tab:dims}, it means that there was
no odd torsion.

Table \ref{tab:hos} gives the Hecke eigenvalue data.  In all cases we
computed $T (\ell ,k)$ for $\ell =2, 3, 5, 7$ and $k=1, 2, 3$ (note
that if $\ell \mid N$ then we actually compute an analogue of the
$U$-operator by choosing different coset representatives). For level
$30$, we also computed operators at $\ell =11, 13$.  We give the
eigenvalue data by presenting the Hecke polynomials of each
eigenclass; that is, we give the polynomial on the left of
\eqref{eqn:hp}.

In all cases we were able to match our eigenclasses to Galois
representations.  In particular, let $\varepsilon$ denote the $p$-adic
cyclotomic character of $G_{\Q}$.  Thus $\varepsilon
(\Frob_{\ell}) = \ell$.  Then in all cases our eigenclasses matched
the Galois representations $1\oplus \varepsilon \oplus
\varepsilon^{2}\oplus \varepsilon^{3}$, except when $N=30$.  For this
level, our data matches that produced by the representation $\beta
\oplus \varepsilon \oplus \varepsilon^{2}\oplus \beta \varepsilon
^{3}$, where $\beta$ is the even character $(\Z /30\Z)^{\times}
\rightarrow (\Z /5\Z)^{\times}$ that takes $2$ to $-1$.  As usual, we
do not know how to prove that these Galois representations are
actually attached to their respective Hecke eigenclasses.

\begin{table}
\begin{center}
\begin{tabular}{|c|c|c||c|c|c||}
\hline 
Level $N$&prime $p$&dimension&Level $N$&prime $p$&dimension\\
\hline
11&5&1&27&3&2\\
19&3&1&29&5&1\\
19&5&1&29&7&1\\
22&5&3&30&5&1\\
23&11&1&31&5&1\\
25&5&2&&&\\
\hline
\end{tabular}
\end{center}
\caption{Odd torsion classes\label{tab:dims}}
\end{table}

\begin{table}
\begin{center}
\begin{tabular}{|p{20pt}|p{150pt}|p{20pt}|p{150pt}|}
\hline
\multicolumn{4}{|l|}{$N=11$, $p=5$}\\
\hline
\hline
$T_{2}$&$1-X^{4}$&$T_{5}$&$1-X$\\
$T_{3}$&$1-X^{4}$&$T_{7}$&$1-X^{4}$\\
\hline
\multicolumn{4}{|l|}{$N=19$, $p=3$}\\
\hline
\hline
$T_{2}$&$1+X^{2}+X^{4}$&$T_{5}$&$1+X^{2}+X^{4}$\\
$T_{3}$&$1-X$&$T_{7}$&$1-X-X^{3}+X^{4}$\\
\hline
\multicolumn{4}{|l|}{$N=19$, $p=5$}\\
\hline
\hline
$T_{2}$&$1-X^{4}$&$T_{5}$&$1-X$\\
$T_{3}$&$1-X^{4}$&$T_{7}$&$1-X^{4}$\\
\hline
\multicolumn{4}{|l|}{$N=22$, $p=5$, eigenclass $1$}\\
\hline
\hline
$U_{2}$&$1-X-X^2-X^{3}-X^{4}$&$T_{5}$&$1-X$\\
$T_{3}$&$1-X^{4}$&$T_{7}$&$1-X^{4}$\\
\hline
\multicolumn{4}{|l|}{$N=22$, $p=5$, eigenclass $2$}\\
\hline
\hline
$U_{2}$&$1+2X+X^{2}-2X^{3}-X^{4}$&$T_{5}$&$1-X$\\
$T_{3}$&$1-X^{4}$&$T_{7}$&$1-X^{4}$\\
\hline
\multicolumn{4}{|l|}{$N=22$, $p=5$, eigenclass $3$}\\
\hline
\hline
$U_{2}$&$1+X-X^{2}+X^{3}-X^{4}$&$T_{5}$&$1-X$\\
$T_{3}$&$1-X^{4}$&$T_{7}$&$1-X^{4}$\\
\hline
\multicolumn{4}{|l|}{$N=23$, $p=11$}\\
\hline
\hline
$T_{2}$&$1-4X+4X^{2}+X^{3}-2X^{4}$&$T_{5}$&$1-2X+4X^{2}+3X^{3}+5X^{4}$\\
$T_{3}$&$1+4X+5X^{2}-2X^{3}+3X^{4}$&$T_{7}$&$1-4X-4X^{2}+3X^{3}+4X^{4}$\\
\hline
\multicolumn{4}{|l|}{$N=25$, $p=5$, both eigenclasses}\\
\hline
\hline
$T_{2}$&$1-X^{4}$&$U_{5}$&$1$\\
$T_{3}$&$1-X^{4}$&$T_{7}$&$1-X^{4}$\\
\hline
\multicolumn{4}{|l|}{$N=27$, $p=3$, both eigenclasses}\\
\hline
\hline
$T_{2}$&$1+X^{4}$&$T_{5}$&$1+X^{2}+X^{4}$\\
$U_{3}$&$1$&$T_{7}$&$1-X-X^{3}+X^{4}$\\
\hline
\multicolumn{4}{|l|}{$N=29$, $p=5$}\\
\hline
\hline
$T_{2}$&$1-X^{4}$&$T_{5}$&$1-X$\\
$T_{3}$&$1-X^{4}$&$T_{7}$&$1-X^{4}$\\
\hline
\multicolumn{4}{|l|}{$N=29$, $p=7$}\\
\hline
\hline
$T_{2}$&$1-X-X^{3}+X^{4}$&$T_{5}$&$1-2X-2X^{2}+2X^{3}+X^{4}$\\
$T_{3}$&$1+2X-2X^{2}-2X^{3}+X^{4}$&$T_{7}$&$1-X$\\
\hline
\multicolumn{4}{|l|}{$N=30$, $p=5$}\\
\hline
\hline
$U_{2}$&$1-2X+X^{3}-X^{4}$&$T_{7}$&$1-2X+2X^{2}-X^{3}-X^{4}$\\
$U_{3}$&$1+2X-2X^{2}-X^{4}$&$T_{11}$&$1+X+X^{2}+X^{3}+X^{4}$\\
$U_{5}$&$1$&$T_{13}$&$1+X-2X^{2}+2X^{3}-X^{4}$\\
\hline
\multicolumn{4}{|l|}{$N=31$, $p=5$}\\
\hline
\hline
$T_{2}$&$1-X^{4}$&$T_{5}$&$1-X$\\
$T_{3}$&$1-X^{4}$&$T_{7}$&$1-X^{4}$\\
\hline
\end{tabular}
\end{center}
\caption{Hecke polynomials for odd torsion classes\label{tab:hos}}
\end{table}

\bibliographystyle{amsalpha_no_mr}
\bibliography{AGM-IV}

\end{document}